\documentclass[preprint,12pt]{elsarticle}
\usepackage{amsmath}
\usepackage{a4wide}
\usepackage[utf8]{inputenc}
\usepackage{amssymb}
\usepackage{amsopn}
\usepackage{epsfig}
\usepackage{amsfonts}
\usepackage{latexsym}
\usepackage{amsthm}
\usepackage{enumerate}
\usepackage[UKenglish]{babel}
\usepackage{verbatim}
\usepackage{color}          %Êýѧ×ÖÌå
\usepackage{mathrsfs}
\usepackage{pgf}
\usepackage{tikz}

\usepackage{mathrsfs}   % This provides an alternative way to typeset script capital letters, using \mathscr{}
\usepackage[all,pdf]{xy}

\usetikzlibrary{arrows,automata}

\DeclareMathAlphabet{\pazocal}{OMS}{zplm}{m}{n}
\newtheorem{theorem}{Theorem}[section]
\newtheorem{lemma}[theorem]{Lemma}
\newtheorem{proposition}[theorem]{Proposition}
\newtheorem{corollary}[theorem]{Corollary}

\theoremstyle{definition}
\newtheorem{definition}[theorem]{Definition}
\newtheorem{example}[theorem]{Example}

\theoremstyle{remark}
\newtheorem{remark}[theorem]{Remark}

\numberwithin{equation}{section}

\newcommand{\R}{\ensuremath{\mathbb{R}}}
\newcommand{\N}{\ensuremath{\mathbb{N}}}
\newcommand{\Z}{\ensuremath{\mathbb{Z}}}

\renewcommand{\t}{{\mathbf{t}}}

\renewcommand{\i}{\ensuremath{\mathbf{i}}}

\newcommand{\I}{\mathcal{I}}

\newcommand{\set}[1]{\left\{#1\right\}}

\newcommand{\Ga}{\Gamma}

\newcommand{\f}{\infty}

\newcommand{\Om}{\Omega}

\begin{document}

\begin{frontmatter}



\title{On the union of homogeneous symmetric Cantor set with its translations}


\author[CU]{Derong~Kong}
\ead{derongkong@126.com}
\address[CU]{College of Mathematics and Statistics, Chongqing University, 401331, Chongqing, P.R.China}

\author[ECNU]{Wenxia~Li}
\ead{wxli@math.ecnu.edu.cn}
%\address[ECNU]{School of Mathematical Sciences, Shanghai Key Laboratory of PMMP, East China Normal University, Shanghai 200062, P.R. China}

\author[ECNU]{Zhiqiang~Wang}
\ead{zhiqiangwzy@163.com}
\address[ECNU]{School of Mathematical Sciences, Shanghai Key Laboratory of PMMP, East China Normal University, Shanghai 200062, P.R. China}

\author[ECUST]{Yuanyuan~Yao\corref{cor1}}
\ead{yaoyuanyuan@ecust.edu.cn} \cortext[cor1]{Corresponding author}
\address[ECUST]{School of Mathematics, East China University of Science and Technology, Shanghai 200237, P.R. China}

\author[NFU]{Yunxiu~Zhang}
\ead{zhyunxiu@163.com}
\address[NFU]{Department of Applied Mathematics, Nanjing Forest University, Nanjing 210037, P.R. China}

\begin{abstract}
Fix a positive integer $N$ and a real number $0< \beta < 1/(N+1)$.
Let $\Gamma$ be the homogeneous symmetric Cantor set generated by the IFS
$$ \Big\{ \phi_i(x)=\beta x + i \frac{1-\beta}{N}: i=0,1,\cdots, N \Big\}. $$
For $m\in\mathbb{Z}_+$ we show that there exist infinitely many translation vectors $\t=(t_0,t_1,\cdots, t_m)$ with $0=t_0<t_1<\cdots<t_m$ such that the union $\bigcup_{j=0}^m(\Ga+t_j)$ is a self-similar set.
Furthermore, for $0< \beta < 1/(2N+1)$, we give a complete characterization on which the union $\bigcup_{j=0}^m(\Ga+t_j)$ is a self-similar set. Our characterization relies on determining whether some related directed graph has no cycles, or whether some related adjacency matrix is nilpotent.
\end{abstract}


\end{frontmatter}


\section{Introduction}\label{sec:introduction}

Self-similar set is a fundamental object in the study of fractal geometry (cf.~\cite{Falconer_1990}).  {A non-empty compact set $E$ in a complete metric space $X$ is called a \emph{self-similar set} if there exists a finite set of contractive similitudes $\mathcal F=\set{f_1, f_2,\cdots, f_n}$ such that  $E=\bigcup_{i=1}^n f_i(E)$}.  The set $\mathcal F$ of contractive similitudes is called an \emph{iterated function system} (simply called, IFS) for the self-similar set $E$ (see \cite{Hutchinson_1981}).
In this paper we will study when the union of a self-similar set with its translations is again a self-similar set.

Fix a positive integer $N$ and a real number $0< \beta < 1/(N+1)$. Let $\Ga=\Ga_{\beta, \set{0,1,\cdots, N}}$ be the self-similar set in $\R$ generated by the IFS
\begin{equation*}
  \label{eq:IFS-phi}
   \Big\{ \phi_i(x)=\beta x + i \frac{1-\beta}{N}: i=0,1,\cdots, N \Big\}.
\end{equation*}
Then $\Ga$ is the unique non-empty compact set satisfying
\[ \Gamma = \bigcup_{i=0}^N \phi_i(\Gamma), \]
and it can be written as
 \begin{equation}\label{eq:Ga}
 \Ga=\set{\frac{1-\beta}{N} \sum_{k=1}^\f j_k \beta^{k-1}: j_k\in\set{0,1,\cdots, N}~\forall k\ge 1 }.
 \end{equation}
 Clearly, $\Ga$ is symmetric, i.e., $\Ga=1-\Ga$.
%The self-similar set $\Ga$ is also called  the \emph{middle $(1-2\beta)$-Cantor set}, which is an affine Cantor set which minimizes the Hausdorff dimension with a given thickness $\frac{\beta}{1-2\beta}$ (cf.~\cite{Kraft-1992, Palis_Takens_1993}).

In the literature there is a great interest in the study of intersections of Cantor set with its translations. When $N=1$, Kraft \cite{Kraft-1994} gave a complete description on when the intersection $\Ga_{\beta,\{0,1\}}\cap(\Ga_{\beta,\{0,1\}}+t)$ is a single point, and Li and Xiao \cite{Li-Xiao-1998} calculated the Hausdorff and packing dimensions of the intersection. In \cite{Deng_He_Wen_2008} Deng, He and Wen studied the self-similarity of the intersection of the middle-third Cantor set with its translation, and gave a necessary and sufficient condition for which the intersection is a self-similar set. This result was later extended by Li, Yao and Zhang in \cite{Li_Yao_Zhang_2011} to the homogeneous symmetric Cantor set $\Ga_{\beta,\set{0,1,\ldots, N}}$ for $0<\beta\le 1/(2N+1)$, and by Kong, Li and Dekking \cite{Kong-Li-Dekking-2010} to $\Ga_{\beta,\set{0,1,\ldots, N}}$ for $1/(2N+1)<\beta<1/(N+1)$.

 On the other hand, we know very little on the self-similarity of the union of $\Ga$ with its translations. More precisely, for a translation vector $\t=(t_0, t_1,\cdots, t_m)\in\R^{m+1}$ we are interested in whether the union
 \[\Ga_{\t}:=\bigcup_{j=0}^m(\Ga+t_j)\] is a self-similar set, where for a set $X$ and $a, b\in\R$ we write $aX+b:=\set{ax+b: x\in X}$.
 Note that the self-similarity is invariant under translations. In other words, if $E\subset\R$ is a self-similar set, then so is its translation $E+t$ for any $t\in\R$. Thus, without loss of generality we  assume throughout the paper that the translation vector $\t=(t_0, t_1, \cdots, t_m)\in\R^{m+1}$ always satisfies $0=t_0<t_1<\cdots<t_m$.
 Note that $\Ga=\bigcup_{i=0}^N \phi_i(\Ga)=\bigcup_{i=0}^N \big( \beta{\Ga}+\phi_i(0) \big)$. Then
\[
\beta^{-n}\Ga=\bigcup_{i_1\cdots i_n\in\set{0,1,\cdots, N}^n}\big( \Ga+\beta^{-n}\phi_{i_1\cdots i_n}(0) \big),
\]
where $\phi_{i_1\cdots i_n}=\phi_{i_1}\circ\cdots\circ\phi_{i_n}$ denotes the composition of maps.
It follows that if $m=N^n-1$ and the translation vector $\t=(t_0,t_1,\cdots,t_m)$ takes the values $\big\{ \beta^{-n}\phi_{\mathbf i}(0): {\mathbf i}\in\set{0,1,\cdots, N}^n \big\}$, then the union $\Ga_\t=\beta^{-n}\Ga$ is a self-similar set. However, for other $m\in\Z_+$ can we find $\t\in\R^{m+1}$ such that $\Ga_\t$ is a self-similar set? Our first result answers this  affirmatively.
\begin{theorem}
  \label{th:existence-m-self-similar}
  Suppose  $0< \beta < 1/(N+1)$. Then for any $m\in\Z_+$ there exist infinitely many translation vectors $\t=(t_0, t_1,\ldots, t_m) \in \R^{m+1}$ with $0=t_0<t_1<\cdots<t_m$ such that $\Ga_\t=\bigcup_{j=0}^m(\Ga+t_j)$ is a self-similar set.
\end{theorem}

Next we consider for which translation vector $\t=(t_0, t_1,\ldots, t_m)\in\R^{m+1}$ the union $\Ga_\t=\bigcup_{j=0}^m(\Ga+t_j)$ is a self-similar set.
Observe that an IFS of $\Ga_\t$ might contain a similitude with negative contraction ratio. This makes our characterization of self-similarity of $\Ga_\t$ more complicated.
To describe the self-similarity of $\Ga_\t$ we first introduce the notation of admissible translation vectors (see Definition \ref{def:admissible-translation} below).

Set
\begin{equation}\label{eq:T}
  T:=\bigcup_{n=1}^\f  \set{ \frac{1-\beta}{N} \sum_{k=1}^{n} j_k \beta^{-k}: j_k \in \{0,1,\cdots, N\} \; \forall 1 \le k \le n }.
\end{equation}
For a translation vector $\t=(t_0, t_1, \cdots, t_m)\in T^{m+1}$ let $\tau_{\t}$ be the smallest integer such that each $t_j, 0\le j\le m$, can be written as
 \begin{equation}\label{eq:t-j}
t_j=\frac{1-\beta}{N}\sum_{k=1}^{\tau_\t}t_{j,k}\beta^{-k} \text{ with } t_{j,k} \in \{0,1,\cdots,N\}.
\end{equation}
Then for $n \ge \tau_\t$ we define
\begin{equation}\label{eq:omega}
  \Omega_{\bf t}^n := \big\{ i_1 \cdots i_n \in \{0,1,\cdots,N\}^n : i_{n+1-k} \le N - s_k \text{ for } 1 \le k \le \tau_\t \big\},
\end{equation}
and its conjugate
\begin{equation}\label{eq:hat-omega}
  \hat{\Omega}_{\bf t}^n := \big\{ i_1 \cdots i_n \in \{0,1,\cdots,N\}^n : i_{n+1-k} \ge s_k \text{ for } 1 \le k \le \tau_\t \big\},
\end{equation}
where $s_k=\max_{0\le j\le m}t_{j,k}$ for $1 \le k \le \tau_\t$.
Clearly,  $i_1 i_2 \cdots i_n \in \Omega_\t^n$ if and only if $(N-i_1)(N-i_2) \cdots (N-i_n) \in \hat\Omega_\t^n$.
Note that for any $\i=i_1 \cdots i_n \in \Omega_\t^n$ and $0\le j\le m$ we have
\begin{align*}
  \phi_{\i}(t_j)=\beta^n t_j+\phi_\i(0)&=\frac{1-\beta}{N}\sum_{k=1}^{\tau_\t}t_{j,k}\beta^{n-k}+\frac{1-\beta}{N}\sum_{k=1}^n i_{n+1-k}\beta^{n-k}\\
  &=\frac{1-\beta}{N}\sum_{k=1}^{\tau_\t}(t_{j,k}+i_{n+1-k})\beta^{n-k}+\frac{1-\beta}{N}\sum_{k=\tau_\t+1}^n i_{n+1-k}\beta^{n-k}\in\Ga.
\end{align*}
Similarly, for any $\i=i_1\ldots i_n\in\hat\Om_\t^n$ and $0\le j\le m$,
\[
\phi_\i(-t_j)=-\beta^n t_j+\phi_\i(0)=\frac{1-\beta}{N}\sum_{k=1}^{\tau_\t}(i_{n+1-k}-t_{j,k})\beta^{n-k}+\frac{1-\beta}{N}\sum_{k=\tau_\t+1}^n i_{n+1-k}\beta^{n-k}\in\Ga.
\]

Let $\mathcal{A}_\t^n$ and $\hat{\mathcal{A}}_\t^n$ be the sets of blocks representing the sets $\set{\phi_\i(t_j): \i\in\Om_\t^n; 0\le j\le m}$ and $\big\{\phi_\i(-t_j): \i\in\hat\Om_\t^n; 0\le j\le m \big\}$, respectively:
\begin{equation}\label{eq:A}
\begin{split}
  \mathcal{A}_{\t}^n &= \{ i_1  \cdots i_{n-\tau_\t}( i_{n+1 -\tau_\t} + t_{j,\tau_\t}) \cdots(i_{n-1} + t_{j,2}) (i_n+ t_{j,1}): \i \in \Omega_{\t}^n, \; 0 \le j \le m \} ,\\
  \hat{\mathcal{A}}_{\t}^n &= \{ i_1  \cdots i_{n-\tau_\t}( i_{n+1 -\tau_\t} - t_{j,\tau_\t}) \cdots(i_{n-1} - t_{j,2}) (i_n- t_{j,1}): {\i} \in \hat\Omega_{\t}^n, \; 0 \le j \le m \} .
\end{split}
\end{equation}
Now we define
\begin{equation}\label{eq:W}
  \mathcal{W}_{\bf t}^n := \mathcal{A}_{\bf t}^n \cup \hat{\mathcal{A}}_{\bf t}^n.
\end{equation}
By the definitions of $\Omega_\t^n$ and $\hat\Omega_\t^n$ it follows that $\Omega_\t^{n} = \{0,1,\cdots, N\}^{n-\tau_\t} \times \Omega_\t^{\tau_\t}$ and $\hat\Omega_\t^{n} = \{0,1,\cdots, N\}^{n-\tau_\t} \times \hat\Omega_\t^{\tau_\t}$.
This implies that
\begin{equation}\label{eq:W-relation}
  \mathcal{W}_{\bf t}^n = \{0,1,\cdots, N\}^{n -\tau_\t} \times \mathcal{W}_{\bf t}^{\tau_\t}\quad\forall n\ge \tau_\t.
\end{equation}

%$ \cup \hat\Omega_\t^n$ the digits $i_1,\cdots, i_{n-\tau_\t}$ can be chosen freely from the digit set $\{0,1,\cdots, N\}$.
%%only the elements at the positions after the $(n+1-\tau_\t)$-th position have restrictions; for other elements we have free choices from the digit set $\{0,1,\cdots, N\}$.
%By the definitions of $\Omega_\t^n$ and $\hat\Omega_\t^n$ it follows that $\Omega_\t^{n} = \{0,1,\cdots, N\}^{n-\tau_\t} \times \Omega_\t^{\tau_\t}$ and $\hat\Omega_\t^{n} = \{0,1,\cdots, N\}^{n-\tau_\t} \times \hat\Omega_\t^{\tau_\t}$ for all $n \ge \tau_\t$.

\begin{definition}
  \label{def:admissible-translation}
  A vector $\t=(t_0, t_1,\cdots, t_m)\in \R^{m+1}$ with $0=t_0<t_1<\cdots <t_m$ is called an \emph{admissible translation vector} if $\t\in T^{m+1}$ and there exists   $\ell\ge \tau_\t$ such that
  \begin{equation*}\label{eq:admissible-symbolic}
  \bigcup_{n=\tau_\t}^\ell \set{0,1,\ldots, N}^{n-\tau_\t}\times\mathcal{W}_\t^{\tau_\t}\times\set{0,1,\ldots, N}^{\ell-n}=\set{0,1,\ldots, N}^\ell.
  \end{equation*}
\end{definition}
According to Definition \ref{def:admissible-translation} it is not easy to verify the admissibility of a translation vector $\t$. In the following we give a more handleable approach by constructing a directed graph. 

Given $\t=(t_0,t_1,\ldots, t_m)\in T^{m+1}$, let $G_\t=(V_\t, E_\t)$ be the directed graph define as follows.
Let $V_\t = \{0,1,\cdots, N\}^{\tau_\t} \setminus \mathcal{W}_\t^{\tau_\t}$.
For two vertices ${\mathbf i}=i_1 i_2 \cdots i_{\tau_\t}, {\mathbf j}=j_1 j_2 \cdots j_{\tau_\t} \in V_\t$, we draw a directed edge   from $\i$ to ${\mathbf j}$  if $i_2 \cdots i_{\tau_\t} = j_1 \cdots j_{\tau_\t -1}$. Then $E_\t$ is the collection of all such directed edges.
%In particular, when $\tau_\t =1$ we have $E_\t = \{({\mathbf i},{\mathbf j}): {\mathbf i},{\mathbf j} \in V_\t \}$.
We say that $G_\t$ has a \emph{cycle} if there exists a directed path  in $G_\t$ starting and ending at the same vertex.
For convenience, we say that the empty graph has no cycles.
For the directed graph $G_\t$ we denote its \emph{adjacency matrix} by $A_\t$. Then $A_\t$ is a $0$-$1$ matrix with the size  $\#V_\t\times \#V_\t$, and an entry $1$ in  $A_\t$ corresponds to a directed edge in $G_\t$. We say that $A_\t$ is \emph{nilpotent} if $A_\t^\ell=0$ for some power $\ell\in\Z_+$.

\begin{proposition}
  \label{prop:chara-admissibility}
  Let $\t=(t_0,t_1,\ldots, t_m)\in T^{m+1}$. The following statements are equivalent.
  \begin{enumerate}
    [{\rm(i)}]
    \item $\t$ is admissible;
    \item $G_\t$ has no cycle;
    \item $A_\t$ is nilpotent.
  \end{enumerate}
\end{proposition}
\begin{proof}
  (i) $\Rightarrow$ (iii). Suppose $A_\t$ is not nilpotent. Then for any $\ell\in\Z_+$ the matrix $A_\t^\ell\ne 0$. This implies that for any $\ell\ge \tau_\t$ there exists a path of length $\ell$ in the directed graph $G_\t$. By the construction of $G_\t$ it follows that for any $\ell\ge \tau_\t$ there exists a word $\i$ of length $\ell$ such that each subword of length $\tau_\t$ in $\i$ belongs to $V_\t=\set{0,1,\ldots, N}^{\tau_\t}\setminus \mathcal{W}_\t^{\tau_\t}$. So,
  \[
  \i\in\set{0,1,\ldots, N}^\ell\setminus \bigcup_{n=\tau_\t}^\ell \set{0,1,\ldots, N}^{n-\tau_\t}\times\mathcal{W}_\t^{\tau_\t}\times\set{0,1,\ldots, N}^{\ell-n},
  \]
  which implies that $\t$ is not admissible by Definition \ref{def:admissible-translation}.

  (iii) $\Rightarrow$ (ii). This follows directly by observing that an entry in $A_\t^\ell$ (say, Row $\i$ and Column $\mathbf j$) corresponds to the number of length $\ell$ paths from vertex $\i$ to $\mathbf j$.

  (ii) $\Rightarrow$ (i). Suppose $G_\t$ has no cycles. If $V_\t=\emptyset$, then $\mathcal{W}_\t^{\tau_\t}=\set{0,1,\ldots, N}^{\tau_\t}$, and by Definition \ref{def:admissible-translation} it is clear that $\t$ is admissible. Now for $V_\t\ne\emptyset$ let $\ell=\tau_\t+\#V_\t$. Arbitrarily take a word $\i=i_1\ldots i_\ell\in\set{0,1,\ldots, N}^\ell$, it suffices to prove that
  \begin{equation}\label{eq:equiv-1}
  i_{n_0}i_{n_0+1}\ldots i_{n_0+\tau_\t-1}\in\mathcal{W}_\t^{\tau_\t}\quad\textrm{for some }1\le n_0\le \ell-\tau_\t+1.
  \end{equation}
    Suppose on the contrary that any block of length $\tau_\t$ in $\i$ does not belong to $\mathcal{W}_\t^{\tau_\t}$. Then $i_ni_{n+1}\ldots i_{n+\tau_\t-1}\in V_\t$ for all $1\le n\le \ell-\tau_\t+1$, and this gives a directed path of length $\ell-\tau_\t=\#V_\t$ in $G_\t$. By  Pigeonhole Principle it follows that $G_\t$ contains a cycle, leading to a contradiction with our assumption.  This proves (\ref{eq:equiv-1}) as desired.
\end{proof}
\begin{remark}
We point out that the characterization of admissibility in Proposition \ref{prop:chara-admissibility} is more handleable.
For example, by using the depth-first search  we can detect the existence of cycles in a directed graph (see \cite{Cormen-Leiserson-Rivest-Stein-2009}).
\end{remark}

For a translation vector $\t=(t_0, t_1,\cdots, t_m)\in\R^{m+1}$ with $0=t_0<t_1<\cdots<t_m$ we define its \emph{conjugate} by $\hat\t=(\hat t_0, \hat t_1,\cdots, \hat t_m)$ where $\hat t_j = t_m - t_{m-j}$ for $0\le j \le m$. Then the elements in the vector $\hat\t$ are also listed in a strictly  increasing order. Furthermore,   $\t$ and $\hat\t$ are conjugate to each other, and by the symmetry of $\Ga$ it follows that
\begin{equation}
  \label{eq:IFS+3}
  (1+t_m)-\Ga_\t =\bigcup_{j=0}^m \big( (1-\Ga)+(t_m-t_j) \big)=\bigcup_{j=0}^m\big( \Ga+\hat t_j \big)=\Ga_{\hat\t}.
\end{equation}
Note that $\Ga_\t$ is a self-similar set if and only if $\Ga_{\hat\t}$ is a self-similar set.
Based on the definition of admissible translation vectors we give a necessary and sufficient condition for the union $\Ga_\t=\bigcup_{j=0}^m(\Ga+t_j)$ to be a self-similar set.
\begin{theorem}
  \label{th:characterization-self-similarity}
 Let $0< \beta < 1/(2N+1)$, and   $\t=(t_0,t_1,\ldots, t_m)\in\R^{m+1}$ with $0=t_0<t_1<\cdots<t_m$. Then $\Ga_\t=\bigcup_{j=0}^m(\Ga+t_j)$ is a self-similar set if and only if either $\t$ or its conjugate $\hat\t$ is an admissible translation vector.
\end{theorem}

As an application, we give an explicit characterization on the self-similarity of $\Ga\cup(\Ga+t)$.
For $x\in\R$ let $\lfloor x \rfloor$  denote its integer part.
\begin{corollary}\label{cor:m=1}
  Let $0< \beta < 1/(2N+1)$ and   $t>0$. Then
  $\Ga\cup(\Ga+t)$ is a self-similar set if and only if \[t = \frac{j(1-\beta)}{N} \beta^{-k}  \] for some $j \in \big\{ 1, 2, \cdots, \lfloor\frac{N+1}{2}\rfloor \big\}$ and $k\in\N$.
\end{corollary}

 Note that a similar result of Corollary \ref{cor:m=1} for $N=1$ and $\beta=1/k$ with $k\in\Z_{\ge 3}$ was obtained in \cite[Theorem 1.1]{Deng-Liu-2011}.

The rest of the paper is arranged as follows. In the next section we give some examples.  In Section \ref{sec:generating IFSs} we describe the generating IFSs of $\Ga_\t$. The proofs of Theorem \ref{th:existence-m-self-similar}  and  \ref{th:characterization-self-similarity} will be given in Section \ref{sec:proofs-thms}.

\section{Examples}\label{sec:examples}

In this section we give some examples to illustrate our main results.

\begin{example}
  Fix a positive integer $N$ and a real number $0< \beta < 1/(2N+1)$.
  For $m \in \Z_+$, let $\t=(t_0,t_1, \cdots, t_m) \in \R^{m+1}$ where $t_0 =0$ and for $1 \le j \le m$,
  \[ t_j = \frac{1-\beta}{N} \sum_{k=1}^{j} \beta^{-k}. \]
  Clearly, $\t \in T^{m+1}$.
  By calculation, we have $\tau_\t = m$, $\Omega_\t^{m} = \{0,1,\cdots,N-1\}^m$, and $\hat\Omega_\t^{m} = \{1,2,\cdots, N\}^m$.
  It follows that \[ \mathcal{A}_\t^{m} = \bigcup_{k=0}^m \{0,1,\cdots,N-1\}^{m-k} \times \{1,2,\cdots, N\}^k, \]
  and \[ \hat{\mathcal{A}}_\t^{m} = \bigcup_{k=0}^m \{1,2,\cdots,N\}^{m-k} \times \{0,1,\cdots, N-1\}^k. \]
  Note that $\mathcal{W}_\t^m = \mathcal{A}_\t^{m} \cup \hat{\mathcal{A}}_\t^{m}$ and $V_\t = \{0,1,\cdots,N\}^m \setminus \mathcal{W}_\t^m$.
  The discussion is split into three cases.

  Case (i): $m=1$. By Corollary \ref{cor:m=1} the set $\Ga_\t$ is a self-similar set.

  Case (ii): $m=2$. It is easy to check that $ \mathcal{W}_\t^m = \{0,1,\cdots,N\}^m$.
  This implies that $G_\t$ is an empty graph.
  By Proposition \ref{prop:chara-admissibility} and Theorem \ref{th:characterization-self-similarity}, we conclude that the set $\Ga_\t$ is a self-similar set.

  Case (iii): $m \ge 3$.
  Let $m' = \lfloor m/2 \rfloor$. If $m$ is odd, we have $(0N)^{m'}0, N(0N)^{m'} \not\in \mathcal{W}_\t^m$, and the cycle $(0N)^{m'}0 \to  N(0N)^{m'} \to (0N)^{m'}0$ is in $G_\t$;
  if $m$ is even, we have $(0N)^{m'}, (N0)^{m'} \not\in \mathcal{W}_\t^m$, and the cycle $(0N)^{m'} \to  (N0)^{m'} \to (0N)^{m'}$ is in $G_\t$.
  Note that the conjugate $\hat\t \in T^{m+1}$. We can check that $\mathcal{W}_{\hat\t}^m=\mathcal{W}_\t^m$ and thus, $G_{\hat\t} = G_\t$ has a cycle.
  By Proposition \ref{prop:chara-admissibility} and Theorem \ref{th:characterization-self-similarity}, we conclude that $\Ga_\t$ is \emph{not} a self-similar set for all $m \ge 3$.
\end{example}

\begin{example}
  Let $N=1$ and $0< \beta < 1/3$.
  Take $\t = (t_0,t_1,t_2,t_3) \in \R^4$ where $t_0=0$, and
  \[ t_1 = (1-\beta)(\beta^{-1} + \beta^{-2}), \; t_2 = (1-\beta)(\beta^{-1} + \beta^{-3}), \; t_3 = (1-\beta)(\beta^{-1} + \beta^{-4}).  \]
  Clearly, $\t \in T^{4}$.
  By calculation, we have $\tau_\t = 4$, $\Omega_\t^{4} = \{0000\}$, and $\hat\Omega_\t^{4} = \{1111\}$.
  It follows that \[ \mathcal{A}_\t^{4} = \{0000,0011,0101,1001\},\; \hat{\mathcal{A}}_\t^{4} = \{ 1111,1100,1010,0110 \}. \]
  Note that $\mathcal{W}_\t^4 = \mathcal{A}_\t^{4} \cup \hat{\mathcal{A}}_\t^{4}$, and
  \[ V_\t =\{0,1\}^4 \setminus\mathcal{W}_\t^4 =\{ 0001,0010,0100,0111,1000,1011,1101,1110 \}. \]
  The directed graph
  \[ G_\t:
  \begin{gathered}
  \xymatrix{
  0001 \ar@/^/[r] & 0010 \ar@/^/[d] \\
  1000 \ar@/^/[u] & 0100 \ar@/^/[l]
  }
  \qquad
  \xymatrix{
  0111 \ar@/^/[r] & 1110 \ar@/^/[d] \\
  1011 \ar@/^/[u] & 1101 \ar@/^/[l]
  }
  \end{gathered}
  \]
  has two cycles.
  Note that $\hat \t \not\in T^4$.
  By Proposition \ref{prop:chara-admissibility} and Theorem \ref{th:characterization-self-similarity}, the set $\Ga_\t$ is \emph{not} a self-similar set.
\end{example}

\section{Generating IFSs of the union $\Ga_\t$}\label{sec:generating IFSs}

Given a self-similar set $E\subset\R$, any IFS $\set{f_i(x)=r_i x+b_i}_{i=1}^n$ with $0<|r_i|<1$ and $b_i\in\R$ satisfying $E=\bigcup_{i=1}^n f_i(E)$ is called a \emph{generating IFS} of $E$ (cf.~\cite{Feng_Wang_2009}). Clearly, a self-similar set has infinitely many generating IFSs.
In this section we describe the generating IFSs of $\Ga_\t$.

\begin{proposition}\label{prop:ratio-generating-IFS}
  Let $0<\beta <1/(2N+1)$, and let $\Ga_\t=\bigcup_{j=0}^m (\Ga+t_j)$ be a self-similar set, where $\t=(t_0, t_1,\cdots, t_m)\in\R^{m+1}$ with $0=t_0<t_1<\cdots<t_m$.
  If $r \Ga_\t + b \subset \Ga_\t$ with $0< |r| <1$, then $|r| = \beta^q$ for some $q \in \Z_+$.
\end{proposition}

Our strategy to prove Proposition \ref{prop:ratio-generating-IFS} is as follows: first we prove that either $\Ga_\t$ or $\Ga_{\hat \t} = 1+ t_m - \Ga$ has a generating IFS which contains a similitude $g(x) = r x$ with $0< r < 1$, see Lemma \ref{lem:IFS-positive-ratio}; next we show that $r = \beta^q$ for some $q \in \Z_+$, and either $\t \in T^{m+1}$ or $\hat\t \in T^{m+1}$, see Lemmas \ref{lem:generating-IFS-r} and \ref{lem:translations}; finally we give a complete characterization of all generating IFSs of $\Ga_\t$, see Lemmas \ref{lem:self-similar-condition}, \ref{lem:g-positive}  and \ref{lem:g-negative}.

\begin{lemma}
  \label{lem:IFS-positive-ratio}
  Let $0 < \beta < 1/(N+1)$, and let $\t=(t_0,t_1,\cdots, t_m)\in\R^{m+1}$ with $0=t_0<t_1<\cdots <t_m$.  If $\Ga_\t$ is a self-similar set, then
  either $\Ga_\t$ or $\Ga_{\hat\t}$ has a generating  IFS containing a similitude $g(x)=r x$ with $0< r < 1$.
\end{lemma}
\begin{proof}
  Suppose that $\set{f_i(x)=r_i x+b_i}_{i=1}^n$ is a generating IFS of $\Ga_\t$. Note that $0\in\Ga\subset\Ga_\t$ and $1+t_m=\max\Ga_\t\in\Ga_\t$. Without loss of generality we assume
  \begin{equation}
    \label{eq:IFS+1}
    0\in f_1(\Ga_\t)=r_1\Ga_\t+b_1\quad\textrm{and}\quad 1+t_m\in f_n(\Ga_\t)=r_n\Ga_\t+b_n.
  \end{equation}
  If $r_1>0$, then  (\ref{eq:IFS+1}) implies $b_1=\min f_1(\Ga_\t)=0$, and thus we are done by taking $g(x)=f_1(x)=r_1 x$. If $r_1<0$, then we consider two cases: $r_n>0$, or $r_n<0$.

  Case (I): $r_1<0$ and $r_n>0$. Then by (\ref{eq:IFS+1}) we have $1+t_m=f_n(1+t_m)=r_n(1+t_m)+b_n$, and thus
  \begin{equation}
    \label{eq:IFS+2}
    (1-r_n)(1+t_m)-b_n=0.
  \end{equation}
 Since $\Ga_\t$ is a self-similar set generated by $\set{f_i(x)=r_i x+b_i}_{i=1}^n$, by (\ref{eq:IFS+3}) it follows that
  \begin{align*}
    \Ga_{\hat\t}&=(1+t_m)-\Ga_\t=\bigcup_{i=1}^n(1+t_m -b_i - r_i\Ga_\t)\\
    &=\bigcup_{i=1}^n \Big( r_i(1+t_m-\Ga_\t)+(1-r_i)(1+t_m)-b_i \Big)\\
    &=\bigcup_{i=1}^n \Big( r_i\Ga_{\hat\t}+(1-r_i)(1+t_m)-b_i \Big).
  \end{align*}
 Then $\Ga_{\hat\t}$ is a self-similar set generated by the IFS \[ \set{\hat f_i(x)=r_i x+(1-r_i)(1+t_m)-b_i }_{i=1}^n. \]
 Note by (\ref{eq:IFS+2}) that $\hat f_n(x)=r_n x+(1-r_n)(1+t_m)-b_n=r_n x$. Then we are done by taking $g(x)=\hat f_n(x)=r_n x$.

 Case (II): $r_1<0$ and $r_n<0$. Then by (\ref{eq:IFS+1}) it follows that $0=f_1(1+t_m)$ and $1+t_m=f_n(0)$. This implies that $f_1\circ f_n(0)=0$, and thus  $f_1\circ f_n(x)=r_1r_n x$ with $r_1r_n>0$.
 Note that $\big\{ f_i \circ f_j (x)= r_i r_j x + r_i b_j + b_i\big\}_{1\le i ,j \le n}$ is also a generating IFS of $\Ga_\t$. Hence, we are done by taking $g(x)=f_1\circ f_n(x)=r_1r_n x$.
 %Since $f_1(\Ga_\t)\subset\Ga_\t$ and $f_n(\Ga_\t)\subset\Ga_\t$, we have $f_1\circ f_n(\Ga_\t)\subset\Ga_\t$. So, we can add $f_1\circ f_n$ to the IFS $\set{f_i(x)=r_i x+b_i}_{i=1}^n$ to form a new IFS, which also generates $\Ga_\t$. Hence, we are done by taking $g(x)=f_1\circ f_n(x)=r_1r_n x$.
\end{proof}

For a finite digit set $D\subset\mathbb Z$ let
\begin{equation}\label{eq:Ga-beta-D}
  \Ga_{\beta,D}:=\set{ \frac{1-\beta}{N} \sum_{k=1}^\f j_k\beta^{k-1}: j_k\in D\; \forall k\ge 1 }.
\end{equation}
Then each $x\in \Ga_{\beta, D}$ can be written as $x=\frac{1-\beta}{N}\sum_{k=1}^\f j_k\beta^{k-1}$ with $j_k\in D$, and the infinite sequence $(j_k)=j_1j_2\cdots \in D^{\Z_+}$ is called a \emph{$D$-coding} of $x$. In general, a point in $\Ga_{\beta, D}$ may have multiple $D$-codings.

%We write $D_1 = \{0,1,\cdots, N\}$, $D_2 = D_1 - D_1 = \{ -N, -(N-1), \cdots, -1 , 0 , 1, \cdots, N\}$, and
%$D_3 = D_1 + D_2 = \{ -N, \cdots, -1 , 0 , 1, \cdots, 2N\}$.

For the rest of this section we always assume that \[ 0 < \beta < \frac{1}{2N+1}. \]
The key in our proof is the following result on unique codings.

\begin{lemma}
  \label{lem:unique-codings}
  Each $x\in\Ga_{\beta, \{0,1,\cdots, N\}}\subset\Ga_{\beta, \{ -N, \cdots, -1 , 0 , 1, \cdots, 2N\}}$ has a unique $\{ -N, \cdots, -1 , 0,$ $1, \cdots, 2N \}$-coding which coincides with its $\{0,1,\cdots, N\}$-coding.
\end{lemma}
\begin{proof}
  Let $x\in\Ga_{\beta, \{0,1,\cdots, N\}}$. Since $0< \beta < 1/(2N+1)$, $x$ has a unique $\{0,1,\cdots, N\}$-coding, say $(i_k)$. Note that $x$ also belongs to $\Ga_{\beta, \{ -N, \cdots, -1 , 0 , 1, \cdots, 2N\}}$. Then $x$ has a $\{ -N, \cdots, -1 , 0 ,$ $1, \cdots, 2N\}$-coding, say $(j_k)$. It suffices to prove that $j_k=i_k$ for all $k\ge 1$.

  For $k \ge 1$, we define $j'_k := -\min\{j_k,0\}$.
  Then $x + \frac{1-\beta}{N} \sum_{k=1}^\f j'_k \beta^{k-1}$ can be written as
  \begin{equation}\label{eq:unique-1}
    \frac{1-\beta}{N} \sum_{k=1}^\f (i_k+j'_k) \beta^{k-1} = \frac{1-\beta}{N} \sum_{k=1}^\f (j_k+j'_k) \beta^{k-1}.
  \end{equation}
  Since $i_k,j'_k \in \{0,1,\cdots, N\}$, we have $i_k + j'_k \in \{0,1,\cdots, 2N\}$ for all $k \ge 1$.
  Note that $j_k + j'_k = j_k$ if $j_k \ge 0$; $j_k + j'_k=0$ if $j_k < 0$.
  Thus, we also have $j_k + j'_k \in \{0,1,\cdots, 2N\}$ for all $k \ge 1$.
  Since $0< \beta < 1/(2N+1)$, each point in $\Ga_{\beta,\{0,1,\cdots,2N\}}$ has a unique $\{0,1,\cdots,2N\}$-coding.
  Then (\ref{eq:unique-1}) implies that $i_k + j'_k = j_k + j'_k$ for all $k \ge 1$.
  So, $j_k = i_k$ for all $k \ge 1$.
\end{proof}

\begin{lemma}\label{lem:reduction}\mbox{}
\begin{enumerate}[{\rm(i)}]
\item If $x+y, 2x+y,\cdots, Nx+y \in \Ga$ for some $y \in \Ga$, then $x \in \Ga_{\beta,\{-1,0,1\}}$.
\item  If $x, 2x,\cdots, Nx \in \Ga$, then $x\in \Ga_{\beta,\{0,1\}}$.
\end{enumerate}
\end{lemma}
\begin{proof}
  (i) Take $y \in \Ga=\Ga_{\beta,\set{0,1,\ldots, N}}$. Then we can write it as
  \begin{equation}\label{eq:y}
     y = \frac{1-\beta}{N} \sum_{k=1}^{\f} y_k \beta^{k-1} \text{ with each } y_k \in \{0,1,\cdots, N\}.
  \end{equation}
  Since $x+y \in \Ga$, we have $x \in \Ga - y \subset \Ga -\Ga = \Ga_{\beta,\{-N,\cdots, -1,0,1,\cdots, N\}}$, which can be written as
  \begin{equation}\label{eq:x} x = \frac{1-\beta}{N} \sum_{k=1}^{\f} x_k \beta^{k-1} \text{ with each } x_k \in \{-N,\cdots, -1,0,1,\cdots, N\}. \end{equation}

  We first have \[ x + y = \frac{1-\beta}{N} \sum_{k=1}^{\f} (x_k+y_k) \beta^{k-1}\in\Ga_{\beta,\set{0,1,\ldots, N}}. \]
  Note by (\ref{eq:y}) and (\ref{eq:x})  that  $(x_k+y_k)_{k=1}^\f$ is a $\{-N,\cdots, -1,0,1,\cdots,2N\}$-coding of $x+y$.
 So, by Lemma \ref{lem:unique-codings} it follows that $x_k + y_k \in \{0,1,\cdots, N\}$ for all $k \ge 1$.

  Next, observe that \[ 2x+y = \frac{1-\beta}{N} \sum_{k=1}^{\f} (2x_k+y_k) \beta^{k-1}\in\Ga_{\beta,\set{0,1,\ldots, N}}. \]
 Since $x_k + y_k \in \{0,1,\cdots, N\}$ for all $k \ge 1$, by (\ref{eq:x}) it follows that $(2x_k + y_k)_{k=1}^\f$ is a $\{-N,\cdots, -1,0,$ $1,\cdots,2N\}$-coding of $2x+y$. Hence,
  by Lemma \ref{lem:unique-codings}  we conclude that $2x_k + y_k \in \{0,1,\cdots, N\}$ for all $k \ge 1$.

  Proceeding this argument $N$ times  we conclude that $N x_k + y_k \in \{0,1,\cdots, N\}$ for all $k \ge 1$.
  Note by (\ref{eq:y}) that $y_k \in \{0,1,\cdots, N\}$.
  Thus, $x_k \in \{-1,0,1\}$ for all $k \ge 1$.
  That is, $x \in \Ga_{\beta,\{-1,0,1\}}$.

  (ii) Taking $y =0$, by (i) we obtain $x \in \Ga_{\beta,\{-1,0,1\}}$. Since $x$ also belongs to $\Ga=\Ga_{\beta, \set{0,1,\ldots, N}}$,
     by Lemma \ref{lem:unique-codings} it follows that $x \in \Ga_{\beta, \set{0,1,\ldots, N}} \cap \Ga_{\beta,\{-1,0,1\}} = \Ga_{\beta,\{0,1\}}$.
\end{proof}

\begin{lemma}
  \label{lem:generating-IFS-r}
  Let $\t=(t_0, t_1,\cdots, t_m)\in\R^{m+1}$ with $0=t_0<t_1<\cdots<t_m$. If $r\Ga_\t\subset\Ga_\t$ with $0< r < 1$, then $r=\beta^q$ for some $q\in\Z_+$.
\end{lemma}
\begin{proof}
  Note by (\ref{eq:Ga}) that $j(1-\beta)\beta^k/N\in\Ga$ for $k\ge 0$ and $j \in \{1,\cdots,N\}$. Take $k$ large enough so that $r(1-\beta)\beta^k<t_1$. Since $r\Ga_\t\subset\Ga_\t=\bigcup_{j=0}^m(\Ga+t_j)$ with $0=t_0<t_1<\cdots<t_m$, it follows that $ j r(1-\beta)\beta^k/N \in\Ga$ for all $j \in \{1,\cdots,N\}$.
  By Lemma \ref{lem:reduction} (ii), we conclude that $r(1-\beta)\beta^k/N \in \Ga_{\beta,\{0,1\}}$.
  Note that $0< r < 1$. By (\ref{eq:Ga-beta-D}) it follows that
  \begin{equation*}
    r=\beta^{q}+\sum_{k=q+1}^\f j_k\beta^k\quad\textrm{with}\quad j_k\in\set{0,1},
  \end{equation*}
where $q\in\Z_+$. So, it suffices to prove that $j_{k}=0$ for all $k\ge q+1$.

Suppose on the contrary there exists $q'>q$ such that
\begin{equation}\label{eq:r1-1}
  r=\beta^q+\beta^{q'}+\sum_{k=q'+1}^\f j_k \beta^k.
\end{equation}
 Note that $(1-\beta)\beta^n(\beta^q+\beta^{q'})\in\Ga$ for any $n\ge 0$. Take $n$ sufficiently large so that  $r(1-\beta)\beta^n(\beta^q+\beta^{q'})<t_1$. Then we obtain
 \begin{equation}\label{eq:r1-2}
   y:=r(1-\beta)\beta^n(\beta^q+\beta^{q'})\in\Ga=\Ga_{\beta,\set{0,1,\cdots,N}}.
 \end{equation}
 On the other hand, by (\ref{eq:r1-1}) it follows that
 \begin{align*}
   y & =(1-\beta)\beta^n(\beta^q+\beta^{q'})\bigg(\beta^q+\beta^{q'}+\sum_{k=q'+1}^\f j_k\beta^k\bigg)\\
   & =\frac{(1-\beta)\beta^n}{N}\bigg( N\beta^{2q}+2N\beta^{q+q'}+N\beta^{2q'}+\sum_{k=q'+1}^\f j_k N \beta^{k+q} + \sum_{k=q'+1}^\f j_k N \beta^{k+q'} \bigg),
 \end{align*}
 which has a $\set{-N,\cdots,-1,0,1,\cdots,2N}$-coding different from its $\set{0,1,\cdots,N}$-coding by (\ref{eq:r1-2}). This leads to a contradiction with Lemma \ref{lem:unique-codings}.
\end{proof}

Recall that
\[ T=\bigcup_{n=1}^\f  \set{ \frac{1-\beta}{N} \sum_{k=1}^{n} j_k \beta^{-k}: j_k \in \{0,1,\cdots, N\} \; \forall 1 \le k \le n }. \]

\begin{lemma}
  \label{lem:translations}
  Let $\Ga_\t=\bigcup_{j=0}^m(\Ga+t_j)$ be a self-similar set, where $\t=(t_0, t_1,\ldots, t_m)\in\R^{m+1}$ with $0=t_0<t_1<\cdots<t_m$. Then either $\t\in T^{m+1}$ or its conjugate $\hat \t\in T^{m+1}$.
  Furthermore, $t_{j+1}-t_j >1$ for all $0\le j <m$.
\end{lemma}
\begin{proof}
  We first assume that $\Ga_\t$ has a generating IFS which contains a similitude $r x$ with $0< r < 1$.
  Then by Lemma \ref{lem:generating-IFS-r} there exists $q\in\Z_+$ such that $\beta^q\Ga_\t\subset\Ga_\t$. Take $u$ sufficiently large so that $\beta^{uq}(1+t_m)<t_1$. Then we obtain
  \begin{equation}
    \label{eq:trans-1}
    \beta^{uq}\Ga_\t\subset\Ga.
  \end{equation}
  Since $t_j\in\Ga+t_j\subset\Ga_\t$ for each $0\le j\le m$, by (\ref{eq:trans-1}) we have $\beta^{uq}t_j\in\Ga$. This together with (\ref{eq:Ga}) implies that
  \begin{equation}
    \label{eq:trans-2}
    t_j=\frac{1-\beta}{N} \beta^{-uq}\sum_{k=1}^\f t_{j,k}\beta^{k-1} \quad\textrm{with}\quad t_{j,k}\in\set{0,1,\cdots,N}.
  \end{equation}
  So, to prove $t_j\in T$ it suffices to prove that $t_{j,k} =0$ for all $k > uq$.

  Suppose on the contrary that $t_{j,k_1}\ne 0$ for some $k_1 > uq$.
  Note that \[ y:=\frac{1-\beta}{N}(N+1-t_{j,k_1}) \beta^{k_1-uq-1} + t_j\in\Ga_\t. \]
  Then by (\ref{eq:trans-1}) we have $\beta^{uq} y \in\Ga=\Ga_{\beta, \set{0,1,\cdots,N}}$. On the other hand, by (\ref{eq:trans-2}) it follows that
  \begin{align*}
    \beta^{uq}y & = \frac{1-\beta}{N}(N+1-t_{j,k_1}) \beta^{k_1-1} + \frac{1-\beta}{N} \sum_{k=1}^\f t_{j,k}\beta^{k-1} \\
    & = \frac{1-\beta}{N}\bigg( \sum_{k\ge 1, k \ne k_1} t_{j,k}\beta^{k-1} + (N+1) \beta^{k_1-1}\bigg),
  \end{align*}
  which has a $\set{-N, \cdots, -1,0,1,\cdots, 2N}$-coding different from its $\set{0,1,\cdots,N}$-coding. This leads to a contradiction with Lemma \ref{lem:unique-codings}. So,
  \begin{equation*}%\label{eq:trans-3}
    t_j=\frac{1-\beta}{N} \sum_{k=1}^{uq} t_{j,k}\beta^{k-uq-1} \in T\quad\textrm{for all }0\le j\le m.
  \end{equation*}

  Take $j\in\set{0,1,\cdots, m-1}$, and let $k_2=\min\set{1\le k\le uq: t_{j+1,k}\ne t_{j,k}}$. Since $t_{j+1}>t_j$, by using $0< \beta < 1/(2N+1)$ it follows that $t_{j+1,k_2} > t_{j,k_2}$.
  Therefore,
  \[
  t_{j+1}-t_j\ge \frac{1-\beta}{N} \beta^{k_2-uq-1} - \frac{1-\beta}{N} \sum_{k=k_2+1}^{uq} N \beta^{k-uq-1} = 1 + \beta^{k_2-uq} \Big( \frac{1-\beta}{N\beta} - 1 \Big) > 1,
  \]
  as desired.

  Next, we assume that $\Ga_{\hat \t}$ has a generating IFS which contains a similitude $r x$ with $0< r < 1$.
  By the above argument, we conclude that $\hat \t \in T^{m+1}$, and $\hat t_{j+1} - \hat t_j > 1$ for all $0 \le j < m$.
  Note that $\hat t_j = t_m -t_{m-j}$ for all $0 \le j \le m$.
  Thus, we also have $t_{j+1} - t_j = \hat t_{m-j} - \hat t_{m-j-1} > 1$ for all $0 \le j < m$.

  By Lemma \ref{lem:IFS-positive-ratio}, either $\Ga_\t$ or $\Ga_{\hat \t}$ has a generating IFS which contains a similitude $r x$ with $0 < r < 1$.
  Thus, we conclude that either $\t\in T^{m+1}$ or its conjugate $\hat \t\in T^{m+1}$.
  In any case, we have $t_{j+1} - t_j > 1$ for all $0 \le j < m$.
\end{proof}

The following lemma states that $\Ga_\t = \bigcup_{j=0}^m (\Ga + t_j)$ is a self-similar set if and only if $\Ga$ can be written as a union of similar copies of $\Ga_\t$.

\begin{lemma}
  \label{lem:self-similar-condition}
  Let $\t=(t_0,t_1,\ldots, t_m)\in\R^{m+1}$ with $0=t_0<t_1<\cdots<t_m$. Then $\Ga_\t$ is a self-similar set if and only if there exists a finite set $\mathcal G$ of similitudes such that
  \begin{equation}
    \label{eq:ssc-1}
    \bigcup_{g\in\mathcal G}g(\Ga_\t)=\Ga.
  \end{equation}
\end{lemma}
\begin{proof}
  The sufficiency is easier, because if (\ref{eq:ssc-1}) holds for some finite set $\mathcal G$, then
  \[\bigcup_{j=0}^m\set{g(x)+t_j: g\in\mathcal G}\] is a generating IFS of $\bigcup_{j=0}^m(\Ga+t_j)=\Ga_\t$.

  For the necessity suppose that $\Ga_\t$ is a self-similar set generated by an IFS $\mathcal F=\set{f_i(x)=r_ix+b_i}_{i=1}^n$.
  By Lemma \ref{lem:translations}, we have
  \[ \delta := \min_{0\le j<m} (t_{j+1}-t_j -1) = \min_{0\le j_1<j_2\le m}\mathrm{dist}(\Ga+t_{j_1}, \Ga+t_{j_2}) >0. \]
  Note that for any $p \ge 1$, \[ \big\{ f_{i_1} \circ f_{i_2} \circ \cdots \circ f_{i_p} (x): 1\le i_1, i_2, \cdots, i_p \le n \big\} \]
  is again a generating IFS of $\Ga_\t$.
  Without loss of generality we may assume that all similarity ratios are sufficiently small so that $|r_i| < \delta/(1+t_m)$ for all $1 \le i \le n$.
  Then we have $\mathrm{diam}\big( f_i(\Ga_\t) \big) = |r_i|\cdot \mathrm{diam}(\Ga_\t)<\delta$.
  Therefore, for each $f\in\mathcal F$ there exists a unique $j\in\set{0,1,\cdots, m}$ such that $f(\Ga_\t)\subset\Ga+t_j$.  Set $\mathcal G:=\{f\in\mathcal F: f(\Ga_\t)\subset\Ga\}$. Then $\bigcup_{g\in\mathcal G}g(\Ga_\t)=\Ga$ as desired.
\end{proof}

Recall that for a translation vector $\t\in T^{m+1}$ the sets $\Omega_\t^n$ and $\hat\Omega_\t^n$ are defined in (\ref{eq:omega}) and (\ref{eq:hat-omega}) for $n \ge \tau_\t$.

\begin{lemma}\label{lem:g-positive}
  Let ${\bf t} =(t_0, t_1, \cdots, t_m) \in T^{m+1}$ with $0=t_0 < t_1 < \cdots < t_m$, and suppose $g(x) = r x + b$ with $0< r < 1$.
  If $g(\Gamma_{\bf t}) \subset \Gamma$, then we have $g(x) = \phi_{\bf i}(x)$ for some ${\bf i} \in \Omega_{\bf t}^n$ where $n \ge \tau_\t$.
\end{lemma}
\begin{proof}
  Note that $r \Gamma_{\bf t} + b \subset \Gamma$ and $0 \in \Gamma_\t$.
  It follows that $b \in \Gamma$, and we write
  \begin{equation}\label{eq:b-1}
    b = \frac{1-\beta}{N}\sum_{k=1}^{\f} b_k \beta^{k-1} \text{ with } b_k \in \{0,1,\cdots,N \}.
  \end{equation}
  Note that $j(1-\beta)/N \in \Ga \subset \Ga_\t$ for all $1\le j \le N$ and $r \Gamma_{\bf t} + b \subset \Gamma$.
  So we have $j r (1-\beta)/N + b \in \Ga$ for all $1 \le j \le N$.
  By Lemma \ref{lem:reduction} (i), we have $r(1-\beta)/N \in \Ga_{\beta,\{-1,0,1\}}$.
  Note that $0< r \le \mathrm{diam}(\Gamma) / \mathrm{diam}(\Gamma_{\bf t}) = 1/(1+t_m) < 1/2$ and $0 < \beta < 1/(2N+1) \le 1/3$.
  By (\ref{eq:Ga-beta-D}) it follows that
  \begin{equation}\label{eq:r-1}
    r = \beta^n + \sum_{k=n+1}^{\f} r_k \beta^k  \text{ with } r_k \in \{ -1, 0 , 1\},
  \end{equation}
  where $n \in \Z_+$.
  In view of (\ref{eq:b-1}) and (\ref{eq:r-1}) we will split our proof into the following three steps:
  (i) we show in (\ref{eq:b-1}) that $b_k =0$ for all $k \ge n+1$, and then $b = \frac{1-\beta}{N}\sum_{k=1}^{n} b_k \beta^{k-1}$;
  (ii) we show in (\ref{eq:r-1}) that $r_k = 0$ for all $k \ge n+1$, and then we have $r = \beta^n$ and thus $g(x) =rx+b= \phi_{b_1 b_2 \cdots b_n}(x)$;
  (iii) we show that $n \ge \tau_\t$ and $b_1 b_2 \cdots b_n \in \Omega_{\t}^n$.

  \textbf{Step 1.} We will show in (\ref{eq:b-1}) that $b_k =0$ for all $k \ge n+1$.
  Fix $k_1 \ge n+1$. Note that $(1-\beta)\beta^{k_1 -n -1} \in \Ga \subset \Ga_\t$ and $r \Ga_\t + b \subset \Ga$.
  It follows that \[ y:=r(1-\beta)\beta^{k_1 -n -1} + b \in \Ga = \Ga_{\beta,\{0,1,\cdots,N\}}. \]
  By (\ref{eq:b-1}) and (\ref{eq:r-1}) we have
  \[ y = \frac{1-\beta}{N} \bigg( \sum_{k=1}^{k_1-1} b_k \beta^{k-1} + (N + b_{k_1}) \beta^{k_1 -1} + \sum_{k=k_1+1}^{\f}(N r_{k+n-k_1} + b_k)\beta^{k-1} \bigg).\]
  By Lemma \ref{lem:unique-codings}, we have $N+b_{k_1} \in \{0,1,\cdots,N\}$.
  This implies $b_{k_1}=0$.
  Since $k_1$ is arbitrary, we conclude that
  \begin{equation}\label{eq:b-2}
    b = \frac{1-\beta}{N}\sum_{k=1}^{n} b_k \beta^{k-1} \text{ with } b_k \in \{0,1,\cdots,N \}.
  \end{equation}

  \textbf{Step 2.} We will show $r=\beta^n$, that is, $r_k = 0$ for all $k \ge n+1$ in (\ref{eq:r-1}).

  First, suppose that there exists $k_2\ge n+1$ such that $r_{k_2} =-1$.
  Note that $1-\beta \in \Ga \subset \Ga_\t$ and $r \Ga_\t + b \subset \Ga$.
  It follows that $r(1-\beta) + b \in \Ga=\Ga_{\beta,\{0,1,\cdots,N\}}$.
  On the other hand, by (\ref{eq:r-1}) and (\ref{eq:b-2}), we have \[ r(1-\beta) + b = \frac{1-\beta}{N}\bigg(\sum_{k=1}^{n} b_k \beta^{k-1} + N \beta^n + \sum_{k=n+1}^{\f} Nr_k \beta^k \bigg), \]
  which has a $\{-N,\cdots,-1,0,1,\cdots,2N\}$-coding different from its $\{0,1,\cdots, N\}$-coding.
  This leads to a contradiction with Lemma \ref{lem:unique-codings}.
  Thus we have $r_k \in \{0,1\}$ for all $k \ge n+1$.

  Next, suppose that there exists $k_3\ge n+1$ such that
  \begin{equation}\label{eq:r-2}
    r = \beta^n + \beta^{k_3} + \sum_{k=k_3+1}^{\f} r_k \beta^k  \text{ with } r_k \in \{ 0 , 1\}.
  \end{equation}
  Note that $(1-\beta)(\beta^n + \beta^{k_3}) \in \Ga \subset \Ga_\t$ and $r \Ga_\t + b \subset \Ga$.
  It follows that $z:=r(1-\beta)(\beta^n + \beta^{k_3}) + b \in \Ga=\Ga_{\beta,\{0,1,\cdots,N\}}$.
  On the other hand, by (\ref{eq:b-2}) and (\ref{eq:r-2}) we have
  \begin{align*}
    z & = (1-\beta)(\beta^n + \beta^{k_3})\bigg( \beta^n + \beta^{k_3} + \sum_{k=k_3+1}^{\f} r_k \beta^k \bigg) +b \\
    & = \frac{1-\beta}{N}\bigg( \sum_{k=1}^{n} b_k \beta^{k-1} + N \beta^{2n} + 2N \beta^{n+k_3} + N \beta^{2k_3} + \sum_{k=k_3+1}^{\f} r_k \beta^{k+n} + \sum_{k=k_3+1}^{\f} r_k \beta^{k+k_3} \bigg),
  \end{align*}
  which has a $\{-N,\cdots,-1,0,1,\cdots,2N\}$-coding different from its $\{0,1,\cdots, N\}$-coding.
  This leads to a contradiction with Lemma \ref{lem:unique-codings}.
  Thus we conclude that $r = \beta^n$.

  \textbf{Step 3.} By the definition of $\tau_\t$, there exists $1\le j_1 \le m$ such that $t_{j_1,\tau_\t} \ne 0$.
  Then we have \[ t_{j_1} = \frac{1-\beta}{N} \sum_{k=1}^{\tau_\t} t_{j_1,k} \beta^{-k} \ge \frac{1-\beta}{N} \beta^{-\tau_\t}. \]
  Note that $t_{j_1} \in \Ga_\t$ and $r\Ga_\t + b \subset \Ga$.
  It follows that $r t_{j_1} + b \in \Ga \subset [0,1]$.
  Thus, $\frac{1-\beta}{N} \beta^{n-\tau_\t} \le r t_{j_1} \le  r t_{j_1} + b \le 1$.
  By using $0<\beta<1/(2N+1)$ this implies that $n \ge \tau_\t$.
  Note that \[ g(x) = \beta^n x + \frac{1-\beta}{N} \sum_{k=1}^{n} b_k \beta^{k-1} = \phi_{b_1 b_2 \cdots b_n}(x). \]
  It remains to show that $b_1 b_2 \cdots b_n \in \Omega_\t^n$, that is, $b_{n+1-k} \le N- s_k$ for all $1 \le k \le \tau_\t$.

  Suppose on the contrary that $b_{n+1-k_4} + s_{k_4} \ge N + 1$ for some $1 \le k_4 \le \tau_\t$.
  By the definition of $s_{k_4}$, there exists $1 \le j_2 \le m$ such that $t_{j_2,k_4} = s_{k_4}$.
  Then we have $b_{n+1-k_4} + t_{j_2,k_4} \ge N + 1$.
  %Let $\i= b_1 b_2 \cdots b_{n-\tau_\t} (b_{n+1-\tau_\t} + t_{j_2,\tau_\t}) (b_{n+2-\tau_\t} + t_{j_2,\tau_\t-1}) \cdots (b_n + t_{j_2, 1})$.
  Note that
  \begin{align*}
    g(\Gamma + t_{j_2})& = \beta^n \Ga + \beta^n t_{j_2} + b \\
     & = \beta^n \Ga + \frac{1-\beta}{N} \bigg( \sum_{k=1}^{n-\tau_\t} b_k \beta^{k-1} + \sum_{k=n+1-\tau_t}^{n}( b_{k} + t_{j_2,n+1-k} ) \beta^{k-1} \bigg).
  \end{align*}
  By Lemma \ref{lem:unique-codings}, we have $g(\Gamma + t_{j_2}) \cap \Ga = \emptyset$.
  This contradicts with $g(\Gamma_{\bf t}) \subset \Gamma$.
  Therefore, we conclude that $b_1 b_2 \cdots b_n \in \Omega_\t^n$, as desired.
\end{proof}

\begin{lemma}\label{lem:g-negative}
  Let ${\t} =(t_0, t_1, \cdots, t_m) \in T^{m+1}$ with $0=t_0 < t_1 < \cdots < t_m$, and suppose $g(x) = -r x + b$ with $0< r < 1$.
  If $g(\Gamma_{\t}) \subset \Gamma$, then we have $g(x) = {\phi_{\i}}(1 - x )$ for some ${\i} \in \hat{\Omega}_{\t}^n$ where $n \ge \tau_\t$.
\end{lemma}
\begin{proof}
  Let $b' = 1- b$.
  Then we have $g(\Ga_\t)= - r \Ga_{\t} + b  = 1- (r \Ga_\t + b') \subset \Ga$.
  Note that $\Ga$ is symmetric.
  It follows that \[ r \Ga_\t + b' \subset 1- \Gamma = \Gamma. \]
  By Lemma \ref{lem:g-positive}, we have \[ r = \beta^n ,\; b' = \phi_{\i'}(0), \]
  where $n \ge \tau_\t$ and ${\i'}=i_1' i_2' \cdots i_n' \in \Omega_{\bf t}^n$.

  Let ${\i} = i_1 i_2 \cdots i_n$ with $i_k = N - i_k'$ for all $1 \le k \le n$.
  Clearly, we have ${\i} \in \hat{\Omega}_{\t}^n$, and
  \begin{align*}
    g(x) & = -\beta^n x + 1 - \phi_{\i'}(0) \\
    & = -\beta^n x + 1 - \frac{1-\beta}{N} \sum_{k=1}^{n} i'_k \beta^{k-1} \\
    & = \beta^n(1-x) + \frac{1-\beta}{N}\sum_{k=1}^{n} (N-i'_k) \beta^{k-1} \\
    & = \beta^n(1-x) + \phi_{\i}(0) \\
    & = \phi_{\i}(1-x),
  \end{align*}
  as desired.
\end{proof}

\begin{proof}[Proof of Proposition \ref{prop:ratio-generating-IFS}]
  Suppose that $\Ga_\t$ is a self-similar.
  By Lemma \ref{lem:translations}, we have either $\t \in T^{m+1}$ or $\hat \t \in T^{m+1}$.
  Furthermore, we have
  \[ \delta := \min_{0\le j<m} (t_{j+1}-t_j -1) = \min_{0\le j_1<j_2\le m}\mathrm{dist}(\Ga+t_{j_1}, \Ga+t_{j_2}) >0. \]
  Write $g(x) = r x + b$.
  Then we have $g^n(\Ga_\t) \subset \Ga_\t$ for all $n \ge 1$.
  Take $n$ large enough so that $|r|^n \cdot \mathrm{diam}(\Ga_\t) < \delta$.
  Then there exist $j_1,j_2 \in \{0,1,\cdots, m\}$ such that
  \begin{equation}\label{eq:g-n-n+1}
    g^n(\Ga_\t) \subset \Ga + t_{j_1}, \; g^{n+1}(\Ga_\t) \subset \Ga + t_{j_2}.
  \end{equation}

  If $\t \in T^{m+1}$, let $g_1(x) = r^n x + g^n(0) - t_{j_1}$ and $g_2(x) = r^{n+1} x + g^{n+1}(0) - t_{j_2}$.
  Then by (\ref{eq:g-n-n+1}) we have $g_1(\Ga_\t) = g^n(\Ga_\t) - t_{j_1} \subset \Ga$ and $g_2(\Ga_\t) = g^{n+1}(\Ga_\t) - t_{j_2} \subset \Ga$.
  By {Lemmas} \ref{lem:g-positive} and \ref{lem:g-negative}, we conclude that $|r|^n = \beta^{q_1}$ and $|r|^{n+1} = \beta^{q_2}$ for some $q_1, q_2 \in \Z_+$.
  It follows that $|r| = \beta ^q$ for some $q \in \Z_+$.

  If $\hat \t \in T^{m+1}$, let $g_3(x) = -r^n x + g^n(1+t_m) - t_{j_1}$ and $g_4(x) = -r^{n+1} x + g^{n+1}(1+t_m) - t_{j_2}$.
  Note that $\Ga_\t = 1+ t_m - \Ga_{\hat \t}$. By (\ref{eq:g-n-n+1}) we have $g_3(\Ga_{\hat\t}) = g^n(1+t_m - \Ga_{\hat\t}) - t_{j_1} = g^n(\Ga_\t) - t_{j_1} \subset \Ga$ and $g_4(\Ga_{\hat\t}) = g^{n+1}(1+t_m - \Ga_{\hat\t}) - t_{j_2}= g^{n+1}(\Ga_\t) - t_{j_2} \subset \Ga$.
  By {Lemmas} \ref{lem:g-positive} and \ref{lem:g-negative}, we conclude that $|r|^n = \beta^{q_3}$ and $|r|^{n+1} = \beta^{q_4}$ for some $q_3, q_4 \in \Z_+$.
  It follows that $|r| = \beta ^q$ for some $q \in \Z_+$.
\end{proof}

\section{Proofs of Theorem \ref{th:existence-m-self-similar} and \ref{th:characterization-self-similarity} }\label{sec:proofs-thms}

In this section we will prove our main theorems. First we prove Theorem \ref{th:characterization-self-similarity}. Recall from Definition \ref{def:admissible-translation} the admissible translation vectors. In the following  we add an equivalent condition for admissible translation vectors to Proposition \ref{prop:chara-admissibility}.

\begin{lemma}\label{lem:admissible-translation-symbolic}
  Suppose that $0< \beta < 1/(N+1)$. A vector ${\bf t}=(t_0, t_1 ,  \cdots , t_m) \in \R^{m+1}$ with $0=t_0 < t_1 < \cdots < t_m$ is an admissible translation vector
  if and only if ${\bf t} \in T^{m+1}$, and there exist finite sets $\I_1 \subset\bigcup_{n\ge \tau_\t}\Om_\t^n$ and $\I_2\subset\bigcup_{n\ge \tau_\t} \hat\Om_\t^n$ such that
   \[ \bigcup_{\i\in\I_1}\phi_\i(\Ga_\t)\cup\bigcup_{\i\in\I_2}\phi_\i(1-\Ga_\t)=\Ga. \]
\end{lemma}
\begin{proof}
Suppose that $\t \in T^{m+1}$. Recall from (\ref{eq:t-j}) that for $0 \le j \le m$, we have \[ t_j = \beta^{-\tau_\t} \phi_{t_{j,\tau_\t} \cdots t_{j,2} t_{j,1}}(0). \]
Take $\i = i_1 i_2 \cdots i_n \in \Omega_{\t}^n$ for $n \ge \tau_\t$.
By the definition of $\Omega_{\t}^n$ in (\ref{eq:omega}), we have
\[ i_1 i_2 \cdots i_{n-\tau_\t}( i_{n+1 -\tau_\t} + t_{j,\tau_\t}) \cdots(i_{n-1} + t_{j,2}) (i_n+ t_{j,1}) \in \{0,1,\cdots,N\}^n \] for all $0\le j \le m$.
It follows that
\begin{align*}
  \phi_{\i}(\Gamma_{\bf t})  = \bigcup_{j =0}^m \phi_{\i}(\Gamma + t_j)
  & = \bigcup_{j =0}^m \big( \beta^{n}\Gamma + \beta^{n} t_j + \phi_{\i}(0) \big) \\
  & = \bigcup_{j =0}^m \big( \beta^{n}\Gamma + \beta^{n-\tau_\t}\phi_{t_{j,\tau_\t} \cdots t_{j,2} t_{j,1}}(0) + \phi_{\bf i}(0) \big) \\
  & = \bigcup_{j =0}^m \big( \beta^{n}\Gamma + \phi_{i_1 i_2 \cdots i_{n-\tau_\t}( i_{n+1 -\tau_\t} + t_{j,\tau_\t}) \cdots(i_{n-1} + t_{j,2}) (i_n+ t_{j,1})}(0) \big) \\
  & = \bigcup_{j =0}^m \phi_{i_1 i_2 \cdots i_{n-\tau_\t}( i_{n+1 -\tau_\t} + t_{j,\tau_\t}) \cdots(i_{n-1} + t_{j,2}) (i_n+ t_{j,1})}(\Gamma).
\end{align*}
Recall the definition of $\mathcal{A}_\t^n$ in (\ref{eq:A}), and so we have
\[\bigcup_{ \i \in \Omega_{{\t}}^n} \phi_\i(\Gamma_{\t}) = \bigcup_{ \i\in \mathcal{A}_{\t}^n } \phi_\i(\Gamma). \]
Similarly, for $\i = i_1 i_2 \cdots i_n \in \hat\Omega_{\t}^n$ where $n \ge \tau_\t$, by the symmetry of $\Ga$ we have
\begin{align*}
  \phi_{\i}(1-\Gamma_{\t}) & = \bigcup_{j =0}^m \phi_{\i}(1-\Gamma - t_j) = \bigcup_{j =0}^m \phi_{\i}(\Gamma - t_j) \\
  & = \bigcup_{j =0}^m \big( \beta^{n}\Gamma - \beta^{n} t_j + \phi_{\i}(0) \big) \\
  & = \bigcup_{j =0}^m \phi_{i_1 i_2 \cdots i_{n-\tau_\t}( i_{n+1 -\tau_\t} - t_{j,\tau_\t}) \cdots(i_{n-1} - t_{j,2}) (i_n- t_{j,1})}(\Gamma).
\end{align*}
Recall the definition of $\hat{\mathcal{A}}_\t^n$ in (\ref{eq:A}), and so we obtain
\[\bigcup_{ \i \in \hat\Omega_{{\t}}^n} \phi_\i(\Gamma_{\t}) = \bigcup_{ \i\in \hat{\mathcal{A}}_{\t}^n } \phi_\i(\Gamma). \]
Note that we set $\mathcal{W}_\t^n = \mathcal{A}_{\t}^n \cup \hat{\mathcal{A}}_{\t}^n$.
Therefore, we conclude that for $n \ge \tau_\t$,
\begin{equation}\label{eq:inclusion}
  \bigcup_{ \i \in \Omega_{\t}^n} \phi_{\i}(\Gamma_{\t}) \cup \bigcup_{ \i \in \hat{\Omega}_{\t}^n } \phi_{\i}(1-\Gamma_{\t}) = \bigcup_{\i \in \mathcal{W}_\t^n } \phi_\i(\Gamma) \subset \Gamma.
\end{equation}

For the necessity, we take $\I_1= \bigcup_{n=\tau_\t}^\ell \Omega_\t^n$ and $\I_2= \bigcup_{n=\tau_\t}^\ell \hat\Omega_\t^n$.
Recall that $\mathcal{W}_\t^n = \{0,1,\cdots, N\}^{n-\tau_\t} \times \mathcal{W}_\t^{\tau_\t}$ in (\ref{eq:W-relation}).
By (\ref{eq:inclusion}), we have
\begin{align*}
  \bigcup_{\i\in\I_1}\phi_\i(\Ga_\t)\cup\bigcup_{\i\in\I_2}\phi_\i(1-\Ga_\t)
  & = \bigcup_{n=\tau_\t}^\ell \bigcup_{\i \in \mathcal{W}_\t^n } \phi_\i(\Gamma) \\
  & = \bigcup_{n=\tau_\t}^\ell\; \bigcup_{\i \in \mathcal{W}_\t^n \times\{0,1,\cdots,N\}^{\ell-n} } \phi_\i(\Gamma) \\
  & = \bigcup_{n=\tau_\t}^\ell \; \bigcup_{\i \in \{0,1,\cdots,N\}^{n-\tau_\t} \times \mathcal{W}_\t^{\tau_\t} \times\{0,1,\cdots,N\}^{\ell-n} } \phi_\i(\Gamma) \\
  & = \bigcup_{\i \in \{0,1,\cdots,N\}^{\ell} } \phi_\i(\Gamma) \\
  & = \Ga,
\end{align*}
{as desired. }

Next, we prove the sufficiency.
We first have $\t \in T^{m+1}$.
Let $\ell$ be the largest length of words in $\I_1 \cup \I_2$.
Clearly, we have $\I_1 \subset \bigcup_{n=\tau_\t}^\ell \Omega_\t^n$ and $\I_2 \subset \bigcup_{n=\tau_\t}^\ell \hat\Omega_\t^n$.
It follows that
\[ \Ga = \bigcup_{\i\in\I_1}\phi_\i(\Ga_\t)\cup\bigcup_{\i\in\I_2}\phi_\i(1-\Ga_\t) \subset \bigcup_{n=\tau_\t}^\ell \Big(\bigcup_{ \i \in \Omega_{\t}^n} \phi_{\i}(\Gamma_{\t}) \cup \bigcup_{ \i \in \hat{\Omega}_{\t}^n } \phi_{\i}(1-\Gamma_{\t}) \Big). \]
Together with (\ref{eq:inclusion}), we have
\begin{equation}\label{eq:equality}
  \Ga = \bigcup_{n=\tau_\t}^\ell \Big( \bigcup_{ \i \in \Omega_{\t}^n} \phi_{\i}(\Gamma_{\t}) \cup \bigcup_{ \i \in \hat{\Omega}_{\t}^n } \phi_{\i}(1-\Gamma_{\t}) \Big) = \bigcup_{n=\tau_\t}^\ell \; \bigcup_{\i \in \{0,1,\cdots,N\}^{n-\tau_\t} \times \mathcal{W}_\t^{\tau_\t} \times\{0,1,\cdots,N\}^{\ell-n} } \phi_\i(\Gamma).
\end{equation}
Since $0< \beta < 1/(N+1)$, the generating IFS of $\Ga$, \[\Big\{ \phi_i(x)=\beta x + i \frac{1-\beta}{N}: i=0,1,\cdots, N \Big\}\] satisfies the strong separation condition.
Therefore, the equality (\ref{eq:equality}) implies
\[ \bigcup_{n = \tau_\t}^\ell \{0,1,\cdots,N\}^{n-\tau_\t} \times \mathcal{W}_{\bf t}^{\tau_\t} \times \{0,1,\cdots,N\}^{\ell-n} = \{0,1,\cdots,N\}^\ell. \]
{Hence, $\t$ is an admissible translation vector}.
We complete the proof.
\end{proof}

First, we prove the sufficiency in Theorem \ref{th:characterization-self-similarity}.
%Recall the definition of admissible translation vectors, see Definition \ref{def:admissible-translation}.

%To describe the self-similarity of $\Ga_\t$ we first introduce the notation of admissible translation vectors.
%Recall the definitions of $T, \Omega_\t^n, \hat\Omega_\t^n$ in (\ref{eq:T}), (\ref{eq:omega}), and (\ref{eq:hat-omega}).

\begin{lemma}\label{lem:sufficiency}
  Suppose that $0< \beta < 1/(N+1)$.
  If $\t=(t_0, t_1,\cdots, t_m)\in \R^{m+1}$ with $0=t_0<t_1<\cdots <t_m$ is an admissible translation vector, then $\Ga_\t = \bigcup_{j=0}^m (\Ga + t_j)$ is a self-similar set.
\end{lemma}
\begin{proof}
   By Lemma \ref{lem:admissible-translation-symbolic}, there exist finite sets $\I_1 \subset\bigcup_{n\ge \tau_\t}\Om_\t^n$ and $\I_2\subset\bigcup_{n\ge \tau_\t} \hat\Om_\t^n$ such that
   \[ \bigcup_{\i\in\I_1}\phi_\i(\Ga_\t)\cup\bigcup_{\i\in\I_2}\phi_\i(1-\Ga_\t)=\Ga. \]
Set $f_\i(x):=\phi_\i(x)$ for $\i\in\I_1$, and set $g_\i(x):=\phi_\i(1-x)$ for $\i\in\I_2$. Then we have
\[ \bigcup_{\i\in\I_1}f_\i(\Ga_\t) \cup \bigcup_{\i\in\I_2}g_\i(\Ga_\t)=\Ga. \]
This implies that $\Ga_\t = \bigcup_{j=0}^m (\Ga + t_j)$ is a self-similar set generated by the IFS
\[ \big\{ f_{\i}(x) + t_j: \;\i \in \I_1, \; 0 \le j \le m \big\} \cup \big\{ g_{\i}(x) + t_j: \;\i \in \I_2, \; 0 \le j \le m \big\}. \]
The proof is completed.
\end{proof}

\begin{proof}[Proof of Theorem \ref{th:characterization-self-similarity}]
  Note that $\Ga_\t$ is a self-similar set if and only if $\Ga_{\hat\t}$ is a self-similar set.
  The sufficiency follows from Lemma \ref{lem:sufficiency}.
  In the following we prove the necessity.

  Suppose that $\Ga_\t$ is a self-similar set.
  Then by Lemma \ref{lem:translations} we have either $\t\in T^{m+1}$ or its conjugate $\hat\t\in T^{m+1}$. Without loss of generality we may assume $\t\in T^{m+1}$.
  By Lemma \ref{lem:self-similar-condition} there exists a finite set $\mathcal G$ of similitudes such that
  \begin{equation}\label{eq:th2-union}
  \bigcup_{g\in\mathcal G}g(\Ga_\t)=\Ga.
  \end{equation}
  By {Lemmas} \ref{lem:g-positive} and \ref{lem:g-negative}, for each $g \in \mathcal{G}$ we have $g(x) = \phi_\i(x)$ for some $\i \in \bigcup_{n=\tau_\t}^\f \Omega_\t^n$, or $g(x) = \phi_\i(1-x)$ for some $\i \in \bigcup_{n=\tau_\t}^\f \hat\Omega_\t^n$.
  Set $\I_1 = \{ \i: g(x) = \phi_\i(x) \in \mathcal{G} \}$ and $\I_2 = \{ \i : g(x) = \phi_\i(1-x) \in \mathcal{G} \}$.
  Then $\I_1$ and $\I_2$ are finite subsets of $\bigcup_{n=\tau_\t}^\f \Omega_\t^n$ and $\bigcup_{n=\tau_\t}^\f \hat\Omega_\t^n$, respectively.
  Furthermore, we have \[ \Ga = \bigcup_{g\in\mathcal G}g(\Ga_\t) = \bigcup_{\i\in\I_1}\phi_\i(\Ga_\t)\cup\bigcup_{\i\in\I_2}\phi_\i(1-\Ga_\t). \]
  {By Lemma \ref{lem:admissible-translation-symbolic}}, $\t=(t_1,t_2,\cdots,t_m)$ is an admissible translation vector.
\end{proof}

Next we prove Theorem \ref{th:existence-m-self-similar}.
\begin{lemma}\label{lem:admissible-translation-reduction}
  Suppose that $0< \beta < 1/(N+1)$, and let $\t=(t_0, t_1 ,  \cdots , t_m) \in T^{m+1}$ with $0=t_0 < t_1 < \cdots < t_m$.
  If $\beta^q \t \in T^{m+1}$ for some $q \in \Z_+$, then $\t$ is an admissible translation vector if and only if $\beta^q \t$ is an admissible translation vector.
\end{lemma}
\begin{proof}
  Write $\t' = \beta^q \t$, that is, $\t = \beta^{-q} \t'$.
  Note that $\t,\t' \in T^{m+1}$.
  We have $\tau_\t = \tau_{\t'} + q$.
  It's easy to check that $\Omega_{\t}^{\tau_t} = \Omega_{\t'}^{\tau_{\t'}} \times \{0,1,\cdots,N\}^q$, and $\hat\Omega_{\t}^{\tau_t} = \hat\Omega_{\t'}^{\tau_{\t'}} \times \{0,1,\cdots,N\}^q$.
  Then we also have $\mathcal{A}_{\t}^{\tau_\t}= \mathcal{A}_{\t'}^{\tau_{\t'}} \times \{0,1,\cdots,N\}^q$ and
  $\hat{\mathcal{A}}_{\t}^{\tau_\t}= \hat{\mathcal{A}}_{\t'}^{\tau_{\t'}} \times \{0,1,\cdots,N\}^q$.
  It follows that \[ \mathcal{W}_\t^{\tau_\t} = \mathcal{W}_{\t'}^{\tau_{\t'}} \times \{0,1,\cdots,N\}^q. \]
  Therefore, for $\ell \ge \tau_\t$ the following three equalities are equivalent:
  \begin{itemize}
    \item $\displaystyle \bigcup_{n = \tau_\t}^\ell \{0,1,\cdots,N\}^{n-\tau_\t} \times \mathcal{W}_{\t}^{\tau_\t} \times \{0,1,\cdots,N\}^{\ell-n} = \{0,1,\cdots,N\}^\ell$;
    \item $\displaystyle \bigcup_{n' = \tau_{\t'}}^{\ell-q} \{0,1,\cdots,N\}^{n'-\tau_{\t'}} \times \mathcal{W}_{\t'}^{\tau_{\t'}} \times \{0,1,\cdots,N\}^{\ell-n'} = \{0,1,\cdots,N\}^\ell$;
    \item $\displaystyle \bigcup_{n' = \tau_{\t'}}^{\ell-q} \{0,1,\cdots,N\}^{n'-\tau_{\t'}} \times \mathcal{W}_{\t'}^{\tau_{\t'}} \times \{0,1,\cdots,N\}^{\ell-q-n'} = \{0,1,\cdots,N\}^{\ell-q}$.
  \end{itemize}
  We conclude that $\t$ is an admissible translation vector if and only if $\t'=\beta^q \t$ is an admissible translation vector.
\end{proof}

\begin{proof}[Proof of Theorem \ref{th:existence-m-self-similar}]
  By Lemma \ref{lem:sufficiency}, it suffices to show that for $m\in \Z_+$ there exist infinitely many admissible translation vectors $\t=(t_0, t_1,\cdots, t_m) \in \R^{m+1}$ with $0=t_0 < t_1 < \cdots < t_m$.
  By Lemma \ref{lem:admissible-translation-reduction}, if $\t \in \R^{m+1}$ is an admissible translation vector, then $\beta^{-k} \t \in \R^{m+1}$ is also an admissible translation vector for all $k \ge 1$.
  Thus, we only need to show that for $m\in \Z_+$ there exists an admissible translation vectors $\t=(t_0, t_1,\cdots, t_m) \in \R^{m+1}$ with $0=t_0 < t_1 < \cdots < t_m$.

  For $m=1$,  take $\t = (0, \frac{1-\beta}{N\beta}) \in T^2$ and then we have $\Ga_\t = \Ga \cup (\Ga + \frac{1-\beta}{N\beta})$.
  It's easy to calculate that $\tau_\t = 1$ and $\Omega_\t^{\tau_\t} = \{0,1,\cdots,N-1\}$.
  Then we have \[ \bigcup_{i \in \Omega_\t^{\tau_\t}} \phi_i(\Ga_\t) = \bigcup_{i \in \Omega_\t^{\tau_\t}} \Big(\beta\Ga + i\frac{1-\beta}{N} \Big) \cup \Big( \beta\Ga + (i+1)\frac{1-\beta}{N} \Big) = \bigcup_{i=0}^N \phi_i(\Ga) = \Ga. \]
  By {Lemma \ref{lem:admissible-translation-symbolic}}, $\t \in T^2$ is an admissible translation vector.

  For $m \ge 2$, there exists $\ell \in \Z_+$ such that $2^\ell \le m < 2^{\ell+1}$.
  Note that $2^\ell +1 \le m+1 \le 2^{\ell+1}$.
  We can find a subset $S \subset \{0,1\}^{\ell+1}$ with $\# S = m+1$ such that $0^{\ell+1}, 1^{\ell+1} \in S$, and for any $i_1 i_2 \cdots i_{\ell+1} \in \{0,1\}^{\ell+1}$ we have either $i_1 i_2 \cdots i_{\ell+1} \in S$ or $(1-i_1) (1-i_2) \cdots (1-i_{\ell+1}) \in S$.
  Let $\t=(t_0, t_1,\cdots,t_m)$ with $0=t_0 < t_1 < \cdots < t_m$ taking the values
  \[ \bigg\{ \frac{1-\beta}{N} \sum_{k=1}^{\ell+1} i_k \beta^{-k}: i_1 i_2 \cdots i_{\ell+1} \in S \bigg\}. \]
  Clearly, we have $\t \in T^{m+1}$.
  Note that $1^{\ell+1} \in S$. We have $\tau_\t = \ell +1$ and $s_k = 1$ for all $1\le k \le \ell+1$.
  It follows that $\Omega_\t^{\ell+1} = \{0,1,\cdots,N-1\}^{\ell+1}$ and $\hat\Omega_\t^{\ell+1} = \{1,2,\cdots,N\}^{\ell+1}$.
  Thus, we have \[ \mathcal{A}_{\t}^{\ell+1} = \big\{ (j_1+i_{\ell+1}) \cdots (j_\ell + i_2) (j_{\ell+1} + i_1):  j_1 j_2 \cdots j_{\ell+1}\in \Omega_\t^{\ell+1}, i_1 i_2 \cdots i_{\ell+1} \in S \big\}, \]
  and
  \[ \hat{\mathcal{A}}_{\t}^{\ell+1} = \big\{ (j_1-i_{\ell+1}) \cdots (j_\ell - i_2) (j_{\ell+1} - i_1):  j_1 j_2 \cdots j_{\ell+1}\in \hat{\Omega}_\t^{\ell+1}, i_1 i_2 \cdots i_{\ell+1} \in S \big\}. \]
  Take $r_1 r_2 \cdots r_{\ell+1} \in \{0,1,\cdots, N\}^{\ell+1}$. Then there exist $j_1 j_2 \cdots j_{\ell+1} \in \{0,1,\cdots,N-1\}^{\ell+1}$ and $i_1 i_2 \cdots i_{\ell+1} \in \{0,1\}^{\ell+1}$ such that
  \[ r_1 r_2 \cdots r_{\ell+1} = (j_1 + i_{\ell+1}) \cdots (j_\ell + i_2) (j_{\ell+1} + i_{1}). \]
  If $i_1 i_2 \cdots i_{\ell+1} \in S$, then we have $r_1 r_2 \cdots r_{\ell+1} \in \mathcal{A}_{\t}^{\ell+1}$;
  if $(1-i_1) (1-i_2) \cdots (1-i_{\ell+1}) \in S$, then $i_1 i_2 \cdots i_{\ell+1} = (1-i'_1) (1-i'_2) \cdots (1-i'_{\ell+1})$ for some $i'_1 i'_2 \cdots i'_{\ell+1} \in S$, and we have
  \[ r_1 r_2 \cdots r_{\ell+1} = (j_1 + 1- i'_{\ell+1}) \cdots (j_\ell + 1 - i'_2) (j_{\ell+1} + 1 - i'_{1}) \in \hat{\mathcal{A}}_{\t}^{\ell+1}. \]
  Thus we conclude that $\mathcal{W}_\t^{\ell+1} = \mathcal{A}_{\t}^{\ell+1} \cup \hat{\mathcal{A}}_{\t}^{\ell+1} = \{0,1,\cdots,N\}^{\ell+1}$.
  This implies that $G_\t$ is an empty graph, which has no cycles.
  By Proposition \ref{prop:chara-admissibility}, $\t \in \R^{m+1}$ is an admissible translation vector.
  We complete the proof.
\end{proof}

Finally we prove Corollary \ref{cor:m=1}.
\begin{proof}[Proof of Corollary \ref{cor:m=1}]
  Let $\t=(t_0,t_1)$ with $0=t_0 < t_1$.
  Note that $\hat\t = \t$.
  By Theorem \ref{th:characterization-self-similarity}, $\Ga_\t = \Ga \cup (\Ga + t_1)$ is a self-similar set if and only if $\t$ is an admissible translation vector.
  We assume $t_1 \in T$ and write \[ t_1 = \frac{1-\beta}{N} \sum_{k=1}^{\tau} s_k \beta^{-k} \text{ with } s_k \in \{0,1,\cdots,N\} \]
  where $\tau \in \Z_+$, and $s_\tau \ne 0$.
  By calculation, we have $\tau_\t = \tau$, and
  \begin{align*}
    \mathcal{W}_\t^\tau = &\; \big\{ i_1 i_2 \cdots i_\tau \in \{0,1,\cdots, N\}^{\tau}: i_{\tau+1-k} \le N -s_k \text{ for } 1 \le k \le \tau \big\} \\
    & \quad \quad \cup \big\{ i_1 i_2 \cdots i_\tau \in \{0,1,\cdots, N\}^{\tau}: i_{\tau+1-k} \ge s_k \text{ for } 1 \le k \le \tau \big\}.
  \end{align*}
  The discussion whether $\t$ is an admissible translation vector is split into three cases.

  \textbf{Case I:} $\tau =1$. Then we have $\mathcal{W}_\t^\tau =\{0,1,\cdots, N-s_1\} \cup \{s_1, s_1 + 1, \cdots, N\}$. In this case, $G_\t$ has no cycles if and only if $\mathcal{W}_\t^\tau = \{0,1,\cdots,N\}$.
  Thus, by Proposition \ref{prop:chara-admissibility}, $\t$ is an admissible translation vector if and only if $s_1 \in \big\{ 1,\cdots, \lfloor\frac{N+1}{2}\rfloor \big\}$.

  \textbf{Case II: } $\tau \ge 2$, and $s_k =0$ for all $1 \le k \le \tau-1$.
  Note that $\beta^{\tau-1} t_1 = \frac{1-\beta}{N} s_\tau \beta^{-1} \in T$.
  By Lemma \ref{lem:admissible-translation-reduction}, $\t$ is an admissible translation vector if and only if $\beta^{\tau-1} \t$ is.
  It follows from \textbf{Case I} that $\t$ is an admissible translation vector if and only if $s_\tau \in \big\{ 1,\cdots, \lfloor\frac{N+1}{2}\rfloor \big\}$.

  \textbf{Case III: }$\tau \ge 2$, and there exists $1\le k \le \tau-1$ such that $s_{k} \ne 0$.
  Let $q=\min\{ 1 \le k \le \tau-1: s_k \ne 0 \}$.
  Note that $$ \beta^{q-1} t_1 = \frac{1-\beta}{N} \sum_{k=q}^{\tau} s_k \beta^{q-k-1} \in T.$$
  By Lemma \ref{lem:admissible-translation-reduction}, $\t$ is an admissible translation vector if and only if $\beta^{q-1} \t$ is.
  Thus we can assume that $s_1 \ne 0$.

  For $i_1 i_2 \cdots i_\tau \in \mathcal{W}_\t^\tau$, we have either $i_1 \le N - s_\tau, i_{\tau}\le N- s_1$, or $i_1 \ge s_\tau, i_{\tau} \ge s_1$.
  It follows that either $i_1, i_\tau \le N - 1$, or $i_1,i_\tau \ge 1$.
  Thus, if $(i_1, i_{\tau}) =(0, N)$, or $(i_1, i_{\tau}) =(N, 0)$, then we have $i_1 i_2 \cdots i_\tau \notin \mathcal{W}_\t^\tau$.
  Note that $V_\t = \{0,1,\cdots, N\}^\tau \setminus \mathcal{W}_\t^\tau$ in the directed graph $G_\t=(V_\t,E_\t)$.
  The following cycle \[ 0^{\tau-1} N \to 0^{\tau-2} N^2 \to \cdots \to 0 N^{\tau-1} \to N^{\tau-1} 0 \to N^{\tau-2} 0^2 \to \cdots \to N 0^{\tau-1} \to 0^{\tau-1} N \] is in $G_\t$.
  By Proposition \ref{prop:chara-admissibility}, $\t$ is not an admissible translation vector.

  Therefore, we conclude that $\Ga \cup (\Ga + t)$ with $t>0$ is a self-similar set if and only if
  \[t = \frac{j(1-\beta)}{N} \beta^{-k}, \] where $j \in \big\{ 1, 2, \cdots, \lfloor\frac{N+1}{2}\rfloor \big\}$ and $k\in\Z_+$.
\end{proof}

\section*{Acknowledgements}

The first author was supported by NSFC No.~11971079. The second author was supported by NSFC No.~12071148 and Science and Technology Commission of Shanghai Municipality (STCSM) No.~18dz2271000.





\end{document}